%% file: main_OR_vRene.tex
\documentclass[nonblindrev]{informs3customR}
\usepackage{appendix}
\usepackage{lmodern}  
\usepackage{amsmath}
\usepackage{mathrsfs}
\usepackage{color}
\usepackage{graphicx}
\usepackage[table]{xcolor} 
\usepackage{setspace}
\usepackage{booktabs}
\usepackage{colortbl}
\usepackage[utf8]{inputenc}
\usepackage[colorlinks=true, linkcolor=blue, citecolor=blue, urlcolor=blue]{hyperref}
\usepackage[capitalise,noabbrev]{cleveref}
\usepackage{mathrsfs}  
\usepackage{natbib}
\bibpunct[, ]{(}{)}{,}{a}{}{,}%
 \def\bibfont{\small}%
\usepackage[english,activeacute]{babel}
\usepackage{graphicx}
\usepackage{amssymb}
\usepackage{setspace}
\usepackage[shortlabels]{enumitem}

\usepackage{bm}

\usepackage{tikz}
\usetikzlibrary{arrows.meta, positioning}
\usepackage{setspace}
\usepackage{float}
\usepackage{pgfplots}
\pgfplotsset{compat=newest} % Sets the compatibility level
\usepackage{mathtools, algorithm, algpseudocode}

\usepackage{enumitem}
\usepackage{bbm}
\usepackage{amsfonts}
\usepackage{comment}
\usepackage{float}

\TheoremsNumberedThrough     
\ECRepeatTheorems
\EquationsNumberedThrough

\makeatletter
\def\footnoterule{\kern-3\p@
  \hrule \@width 6.5in \kern 2.6\p@} % the \hrule is .4pt high
\makeatother

\definecolor{orange}{rgb}{1,0.5,0}
\definecolor{OliveGreen}{cmyk}{0.64,0,0.95,0.40}
\definecolor{BrickRed}{cmyk}{0,0.89,0.94,0.28}

\newcommand{\yichen}[1]{\textcolor{blue}{[Yichen: #1]}}

\newcommand{\ti}[1]{\scalebox{0.45}{\mbox{\rm #1}}} % Adjust the scale as needed
 
\newcommand{\etiny}{\fontsize{5pt}{5pt}\selectfont}
\newcommand{\ttti}[1]{\scalebox{0.7}{\mbox{\rm #1}}}

\newcommand{\tpsi}{\tilde{\psi}}
\newcommand{\hpsi}{\widehat{\psi}}
\newcommand{\eps}{\epsilon}
\newcommand{\var}{\mathbb{V}\mbox{\rm ar}}

\newcommand{\D}{\mbox{\rm d}} 

\newcommand{\kk }{k}
\newcommand{\hpsim}{\hpsi^m}

\newcommand{\GD}{\mathcal{GD}}

\newcommand{\Su }{\mbox{\rm \tiny S}}

\newcommand{\I}{\mbox{\rm \tiny I}}
\newcommand{\FI}{\mbox{\rm \tiny FI}}

\newcommand{\R}{\mbox{\rm \tiny R}}

\newcommand{\e}{\mathbb{E}}

\newtheorem{cor}{Corollary}

\newtheorem{thm}{Theorem}
\newtheorem{exm}{Example}

\renewcommand{\thefootnote}{\fnsymbol{footnote}}
\begin{document}
\RUNAUTHOR{}
\RUNTITLE{}
%\TITLE{\bf Optimizing Information Delays \\ in Supply Chains}
\TITLE{\bf \vspace{-2.7cm}\\Designing Information Delays in Supply Chains \\ }
\ARTICLEAUTHORS{%
\AUTHOR{Prem Talwai}
\AFF{Operations Research Center, Massachusetts Institute of Technology, Cambridge, MA 02142,
\EMAIL{talwai@mit.edu}}
\AUTHOR{Ren\'e Caldentey}
\AFF{Booth School of Business, The University of Chicago, Chicago, IL 60637,
\EMAIL{rene.caldentey@chicagobooth.edu}}
\AUTHOR{Avi Giloni}
\AFF{Sy Syms School of Business, Yeshiva University, New York, NY 10033,
\EMAIL{agiloni@yu.edu}}
\AUTHOR{Clifford Hurvich}
\AFF{Stern School of Business, New York University,  New York, NY 10012,
\EMAIL{churvich@stern.nyu.edu}}
\AUTHOR{David Simchi-Levi}
\AFF{Institute of Data, Systems, and Society, Massachusetts Institute of Technology, Cambridge, MA 02142,
\EMAIL{dslevi@mit.edu}}
\AUTHOR{Yichen Zhang}
\AFF{Daniels School of Business, Purdue University,  West Lafayette, IN 47907,
\EMAIL{zhang@purdue.edu}}
% Enter all authors
}

\date{}
\ABSTRACT{ \vspace{0.5cm}

This paper studies how a downstream retailer in a decentralized two-tier supply chain can implicitly transmit demand information to an upstream supplier through the structure of its order stream in the absence of an explicit information-sharing mechanism. We distinguish our work from prior work by introducing the notion of information delay and by linking optimal implicit information sharing to the group delay of the retailer’s ordering transfer function. We show that pure delay is strictly suboptimal, while fractional-delay mechanisms can reshape the order autocorrelation to improve supplier forecastability and reduce system-wide inventory costs. Using Hardy-space factorization, we develop a tractable family of invertible ARMA policies that approximates the theoretically optimal (but non-rational) limiting filter derived by \cite{caldentey2024information} and preserves its informational delay properties. This construction yields sharp guidance on how policy complexity, as measured by the degrees of the ARMA policies, impacts supply chain costs. We further extend the analysis to memory-constrained suppliers and characterize how the complexity of the retailer’s policy should scale with the supplier’s finite forecasting window, highlighting when, perhaps counterintuitively, increasing policy complexity can become counterproductive.
}
\vspace{0.5cm}

\maketitle

\vspace{-2.1cm} 

% !TEX root = /main_OR.tex 

\section{Introduction}\label{sec:Introduction}

Successful supply chain operations depend on the effective coordination between inventory management strategies and the flow of information that informs them. Synchronizing these physical and informational flows ensures that goods move efficiently through the supply chain, minimizing costs, reducing lead times, and maintaining optimal inventory levels. While inventory policies dictate the timing and quantity of replenishment across firms, their effectiveness is shaped by the availability and accuracy of demand information. When downstream demand is known in a timely and accurate manner, upstream suppliers can anticipate needs, allocate resources efficiently, and maintain service levels. In ideal settings, this information is explicitly shared through collaborative systems or formal communication protocols. In practice, however, such sharing is often absent, delayed, or incomplete—leaving upstream firms to infer demand conditions indirectly from observed order patterns.

This paper investigates how inventory control decisions, and specifically the structure of the retailer’s ordering policy, can serve as a channel for implicit information transmission in the absence of an explicit information-sharing mechanism. We examine how a retailer can design its order stream to convey relevant information to an upstream supplier, enabling improved forecasting and inventory planning.

We build on the results of \cite{caldentey2024information}, who formalized a mathematical framework for evaluating the costs associated with information asymmetries in a simple two-tier supply chain. That work showed that when demand information is fully shared, the retailer can implement a straightforward and intuitive ordering policy that minimizes system-wide inventory costs. In contrast, under no information sharing, the optimal policy becomes significantly more complex and less transparent, involving infinite-memory strategies and impractical non-ARMA signal structures that embed demand information in a way that is only partially recoverable.

Motivated by the limitations of the no-information-sharing regime studied in \cite{caldentey2024information}, this paper revisits the fundamental question of how information can be most effectively transmitted from downstream to upstream stages in supply chains and explores ways to sustain coordinated decision-making when explicit communication is absent. We propose and formalize the notion of \emph{information delay}—a mechanism through which the retailer intentionally slows its replenishment response to demand shocks to transparently embed demand-related information into its order stream. By doing so, the retailer enables the supplier to better infer future demand and make more informed inventory decisions. This implicit signaling mechanism can serve as a substitute for direct information sharing, facilitating alignment in decentralized supply chains.

Our first contribution is to demonstrate that the design of such delays—both in form and intensity—has significant implications for supply chain performance. Departing from prior work (see \citealp{bensoussan2006optimality}) that models delay using a single integer lag (i.e., pure delay), we show, perhaps counterintuitively, that pure delay does not help the supplier forecast perfectly, and that more flexible and effective delay structures can be achieved through the simultaneous use of multiple lags (i.e., smoothing \cite{disney2004variance, BGP2004, brown2004smoothing, caldentey2026information}). At the heart of our analysis is the concept of \emph{fractional delay}, which captures the virtual timescale over which demand information is gradually deferred. We operationalize this through a class of filtering strategies that spread information across time, making order patterns more predictable. This effect is quantified using the concept of \emph{group delay}, which measures how signals are temporally embedded in the order stream. Crucially, we find that unlike pure delay—which leaves the statistical profile of the order process unchanged and is therefore undetectable to the supplier—fractional delay reshapes the autocorrelation structure of orders and improves forecastability. As a result, we establish that in settings without information sharing, pure delay is strictly suboptimal for minimizing joint inventory costs.

Second, we show that while some of the theoretically optimal delay mechanisms are non-invertible and thus impossible for the supplier to detect under no information sharing, they can be closely approximated by practically implementable invertible (i.e., detectable) filters. These approximations yield simple policies that are nearly optimal in terms of minimizing joint inventory costs, while being practical to deploy in decentralized environments.

Third, we extend the analysis to a more realistic setting in which the supplier is memory-constrained and has access only to a finite recent history of past orders.
 This limited visibility introduces an additional layer of friction, diminishing the supplier’s ability to reconstruct the embedded information delay and complicating the use of delay-based signaling. In this context, we uncover and analytically characterize a tradeoff between the complexity of the retailer’s signal encoding and the supplier’s capacity to interpret it. This tradeoff resembles a bias-variance tension between supplier and retailer costs and can be optimized by carefully calibrating the complexity of the delay mechanism to the supplier’s observation horizon. Doing so minimizes overall system costs by balancing precision and tractability in the learning process.

Finally, our paper contributes methodologically to the literature on inventory and supply chain management by reframing information sharing  through the lens of delay embedded in the retailer’s order stream. While prior work \citep{bensoussan2006optimality, chen1999decentralized} focused on pure delay, we study a broader class of all-pass filters that preserve nondeterminism in orders and allow the retailer to shape how demand information is distributed over time. This perspective yields a crisp characterization of the supplier--retailer cost tradeoff in terms of the \emph{type} and \emph{amount} of delay induced by the retailer’s ordering filter. We further show that a crucial all-pass filter can be effectively approximated by a family of minimum-phase (invertible) ARMA filters (Lemma~\ref{Singular Inner Approximation}) that replicate the desired delay properties even when explicit information sharing is absent. Although approximations of all-pass filters by maximum-phase filters have been considered before \citep{mehta2023all, appaiah2020all}, to the best of our knowledge this is the first work to develop minimum-phase (invertible) approximations for this supply-chain signaling problem.

In summary, this paper develops a unifying perspective on information delay, offering a new operational interpretation of demand signaling in supply chains. We show that even without direct communication, the retailer can effectively serve as an information designer, embedding the structure of the demand process into its order stream in a form that is both recoverable and useful to the supplier. By doing so, we connect insights from forecasting theory, signal processing, and inventory management, and introduce new tools for enhancing coordination in decentralized settings. Our findings ultimately suggest that replenishment policies should be carefully designed not only to manage inventory, but also to convey demand information.

\setcounter{footnote}{1}
\section{ A Two-Tier Supply Chain Model} \label{sec:ProblemFormulation}

We consider a prototypical two-echelon supply chain, consisting of a single retailer and a single supplier, evolving in discrete time, as illustrated in \cref{fig:twotierSC}. 
\input{Fig_System_OR}
The operating characteristics mirror those analyzed in \cite{Gavirneni1999,LST2000,Aviv2001,Li-Lee2009,GHS,caldentey2024information}, including: (i) market demand at the retailer is fully backlogged, so all unmet demand is eventually fulfilled; (ii) the supplier follows a state-dependent base-stock policy; and (iii) when facing a shortfall, the supplier can expedite production at an additional cost to maintain a 100\% service level to the retailer.  

The retailer faces uncertain customer demand $D_t$, modeled as an i.i.d.\ Gaussian sequence $D_t = d + \epsilon_t$, where $d>0$ is the mean demand per period and $\{\epsilon_t\}$ is Gaussian white noise with variance $\var(\epsilon_t) = \sigma^2_\epsilon$. The retailer serves demand using its net inventory $I^{\R}_t$ and replenishes by placing periodic orders $O_t$ with the supplier. We restrict attention to weakly stationary replenishment policies with respect to $\{\epsilon_t\}$, so the retailer’s orders admit a one-sided Gaussian MA($\infty$) representation
\begin{equation}
\label{eq: MA Rep}
O_t = d + \sum_{n = 0}^{\infty} \tpsi_n \eps_{t - n},
\end{equation}
where the MA coefficients $\{\tpsi_n\}$ encode the retailer's inventory policy. We write $\tpsi=\{\tpsi_n\}_{n\geq 0}$ for the retailer's policy and note that this class includes most standard replenishment rules, such as base-stock and myopic order-up-to policies.

The supplier fulfills the retailer’s orders using its own inventory $I^{\Su}_t$, which it manages with a state-dependent base-stock policy. In each period the supplier sets its base-stock level $S_t$ equal to the one-step-ahead forecast of next period’s order plus a safety stock:\footnote{Specifically, the supplier determines $S_t$ by solving  
$\displaystyle
\min_S C^{\ti S}=\mathbb{E} \big[h^{\ti S} (S_t - O_{t+1})^+ + b^{\ti S} (S_t - O_{t+1})^- \mid {\cal F}^{\ti S}_t \big].$}
\begin{equation}\label{eq:basestock}
S_t = m_t(\tpsi) + \sigma_{\Su}(\tpsi)\,\zeta^{\Su}, 
\qquad 
\zeta^{\Su} = \Phi^{-1}\!\left(\frac{b^{\Su}}{h^{\Su}+b^{\Su}}\right),
\end{equation}
where $m_t(\tpsi)=\e[O_{t+1}\mid \mathcal{F}^{\Su}_t]$ is the supplier’s one-step-ahead forecast, $\sigma^2_{\Su}(\tpsi)=\var[O_{t+1}\mid \mathcal{F}^{\Su}_t]$ its mean squared forecast error (MSFE), $(h^{\Su},b^{\Su})$ the per-unit holding and expediting costs, and \mbox{$\mathcal{F}^{\Su}_t := \sigma\big(O_\tau : \tau \le t\big)$} denotes the supplier’s information set generated by the order history. The safety factor $\zeta^{\Su}$ balances holding and expediting, and $\Phi$ denotes the standard normal CDF.

Under these operating conditions, the retailer and supplier’s net inventories evolve as
\[
I^{\R}_t = I^{\R}_{t-1} + O_{t-1} - D_t 
\qquad \mbox{and}\qquad
I^{\Su}_t = \big(S_{t-1} - O_t\big)^+.
\]
\cite{caldentey2024information} show that, for any given policy $\tpsi$, the long-run average costs of the retailer and supplier can be written as ${\cal C}^{\R}(\tpsi)=K^{\R}\,\sigma_{\I}(\tpsi)$ and ${\cal C}^{\Su}(\tpsi)=K^{\Su}\,\sigma_{\Su}(\tpsi)$, where $\sigma_{\I}(\tpsi)$ is the stationary standard deviation of the retailer’s inventory, $\sigma_{\Su}(\tpsi)$ is the supplier’s root MSFE, and $K^{\R},K^{\Su}>0$ depend only on cost parameters and safety-stock factors.\footnote{For firm $j\in\{\rm R,S\}$, $K^{j}= \big[ h^{j}\,\zeta^{j}+ (h^{j}+b^{j} )\,{\cal L}(\zeta^{j})\big]$, where $\zeta^{j} = \Phi^{-1}\!\big(\frac{b^{j}}{h^{j}+b^{j}}\big)$ and ${\cal L}(x)=\phi(x)-x(1-\Phi(x))$ is the standard normal loss function, with $\phi$ and $\Phi$ denoting the standard normal PDF and CDF.} Thus, the system-wide cost is
\[
{\cal C}(\tpsi)
={\cal C}^{\R}(\tpsi)+{\cal C}^{\Su}(\tpsi)
=K^{\R}\,\sigma_{\I}(\tpsi)+K^{\Su}\,\sigma_{\Su}(\tpsi).
\]
To ensure stability of the retailer’s inventory process (and hence finiteness of ${\cal C}^{\R}(\tpsi)$), we restrict attention to policies $\tpsi$ that belong to the set of \emph{admissible} policies
\[
\widetilde{\Psi} = \left\{ \tpsi \in \ell^2 : \sum_{n=0}^\infty \tilde{\psi}_n = 1,\quad \sum_{n \geq 0} \left( \sum_{k = 0}^{n-1} \tilde{\psi}_k - 1 \right)^2 < \infty \right\}.
\]
The condition $\sum_{n=0}^\infty \tilde{\psi}_n = 1$ is the \emph{bounded-inventory condition}.

Following \cite{caldentey2024information}, we find it convenient to reformulate the problem of finding an admissible ordering policy $\tpsi \in \widetilde{\Psi}$ that minimizes the supply chain cost ${\cal C}(\tpsi)$ from the time domain into the $z$-transform domain. In particular, we express ${\cal C}(\tpsi)$ in terms of the transfer function $\tpsi(z)$ rather than its time-domain MA coefficients $\{\tpsi_n\}_{n\geq 0}$, where
\[
\tpsi(z) = \sum_{n=0}^\infty \tpsi_n z^n, \qquad z\in\mathbb{D},
\]
and $\mathbb{D}=\{z\in\mathbb{C}:|z|<1\}$ is the open unit disk. Crucially, the admissibility condition $\tpsi\in\widetilde{\Psi}$ implies $\tpsi\in\mathbb{H}^2$, the Hardy--Hilbert space of analytic functions on $\mathbb{D}$ with square-integrable boundary values on the unit circle, which will allow us to exploit a rich set of results from $\mathbb{H}^2$ to characterize an optimal policy. From the inventory dynamics $I^{\R}_t = I^{\R}_{t-1} + O_{t-1} - D_t$ and the fact that the $z$-transform of the demand shock process is the constant 1, the $z$-transform of the retailer’s inventory process satisfies
\[
\psi_{\I}(z) = \frac{z\,\tpsi(z) - 1}{1 - z}.
\]
It follows that the stationary inventory variance can then be written as the energy of $\psi_{\I}$:
\begin{equation}\label{eq_SigmaI_z}
\sigma^2_{\I}(\tpsi) 
= \frac{1}{2\pi} \int_{-\pi}^{\pi} 
\left| \frac{z\,\tpsi(z) - 1}{1 - z} \right|^2 \,\D\lambda,
\qquad z = \exp(-i\lambda).
\end{equation}
Similarly, letting $f_{\tpsi}(\lambda)$ denote the spectral density of the order process induced by $\tpsi$, Kolmogorov’s formula \citep[Section 5.8]{BrockwellDavis} yields
\begin{equation}\label{eq:Kolmogorov}
    \sigma^2_{\Su }(\tpsi) 
    = 2\pi\, \exp\left(\frac{1}{2\pi} \int_{-\pi}^{\pi} \log f_{\tpsi}(\lambda)\, \mathrm{d}\lambda \right),
\end{equation}
and, for our MA($\infty$) representation with i.i.d.\ shocks, 
$f_{\tpsi}(\lambda)=\frac{\sigma^2_\epsilon}{2\pi}\big|\tpsi(e^{-i\lambda})\big|^2$.

For the subsequent analysis, and without loss of generality, we normalize by setting $\sigma_\eps = 1$ and $K^{\Su} = 1$, and defining $K^{\R}=\kappa$. Equivalently, we divide ${\cal C}(\tpsi)$ by $K^{\Su}\sigma_\eps$ and set $\kappa := K^{\R} / K^{\Su}$. This preserves the structure of the problem and makes $\kappa$ the relative weight of the retailer’s cost. Under this normalization, using \eqref{eq_SigmaI_z} and \eqref{eq:Kolmogorov}, the problem of minimizing the total supply chain cost over admissible policies becomes
\begin{equation}\label{eq:Cost_SC_z}
    {\cal C}^*
    = \inf_{\tpsi \in \widetilde{\Psi}}\; 
    {\cal C}(\tpsi) 
    = \inf_{\tpsi \in \widetilde{\Psi}}\;
    \kappa\,\left(\frac{1}{2\pi} \int_{-\pi}^{\pi} 
        \left| \frac{z\,\tpsi(z) - 1}{1 - z} \right|^2\,\D\lambda\right)^{1/2}
    + 
    \exp\!\left(
        \frac{1}{2\pi} 
        \int_{-\pi}^{\pi} 
        \log\big|\tpsi(z)\big|\, 
        \mathrm{d}\lambda 
    \right), \quad z=e^{-i\lambda}.
\end{equation}
The retailer’s policy thus plays a dual role: it shapes both the retailer’s inventory variability and the supplier’s forecast accuracy, and hence the overall supply chain performance.

\section{Prior Results and Overview of Contributions}\label{sec:preliminaries}

We next summarize the solution to \eqref{eq:Cost_SC_z} developed in \cite{caldentey2024information}, which exploits Hardy-space methods, most notably the Smirnov--Beurling (Nevanlinna) inner--outer factorization theorem (see \citealp[Theorem 17.17]{rudin1970real} or \citealp[Theorem 2.6.5]{Nikolski}).
 Any nonzero $\tpsi\in\mathbb{H}^2$ admits a unique (up to a unimodular constant) decomposition $\tpsi(z)=Q(z)\,I(z)$, where $Q$ is \emph{outer} and $I$ is \emph{inner}. The inner factor satisfies $|I(e^{-i\lambda})|=1$ for a.e.\ $\lambda\in[-\pi,\pi)$, and can itself be factored as $I(z)=B(z)\,S(z)$, with $B$ a Blaschke product collecting the zeros of $\tpsi$ inside $\mathbb{D}$ and $S$ a singular inner function determined by a singular measure on the unit circle. Inner functions act as all-pass filters on $\mathbb{T}:= \{z \in \mathbb{C} : |z| = 1\}$, the unit circle in the complex plane; they preserve magnitude while altering phase.

For the supplier’s forecasting problem, this factorization has a clear interpretation. Because inner factors are all-pass, the spectral densities associated with $\tpsi$ and its outer factor $Q$ coincide, so $f_{\tpsi}(\lambda)=f_Q(\lambda)$ almost everywhere in $[-\pi,\pi)$ and, in turn,
\[
\sigma^2_{\Su}(\tpsi)=\sigma^2_{\Su}(Q).
\]
Thus, the outer factor $Q$ captures the \emph{forecastable} (visible) component of the retailer’s orders, whereas the inner factor encodes an \emph{unforecastable} (invisible) component from the supplier’s perspective. In the time domain this corresponds to invertibility: the orders are invertible (the supplier can recover the shocks $\eps_t$ from $\{O_n:n\le t\}$) if and only if $\tpsi$ has no inner factor, i.e., $\tpsi(z)=Q(z)$. %\prem{Maybe add pure delay as example of invisible impact on forecasting costs}

Building on this structure, \cite{caldentey2024information} solve \eqref{eq:Cost_SC_z} by exploiting the interplay between outer and singular inner functions and constructing a sequence of $\epsilon$-optimal policies. Their main result can be summarized as follows.

\begin{thm}\label{thm_CGHZ}
Consider the problem of minimizing ${\cal C}(\tpsi)$ in \eqref{eq:Cost_SC_z} over $\widetilde{\Psi}$. Two regimes arise depending on the retailer’s relative cost weight $\kappa$:
\begin{enumerate}
\item If $\kappa \geq \sqrt{5}$, the optimal policy $\tpsi^*$ is MA(1) and admits the representation
\[
\tpsi^*(z) = \tpsi^*_0 + (1 - \tpsi^*_0)\,z,
\qquad 
\tpsi^*_0 = 1 - \big(\sqrt{\kappa^2 - 1}\big)^{-1}.
\]
This solution is invertible, i.e., has no inner factor.

\item If $\kappa < \sqrt{5}$, let $\gamma_\kappa \geq 0$ be the unique solution of $\kappa^2=(5+2\,\gamma)\,\exp(-2\,\gamma)$ and define, for each $k \in \mathbb{N}$, the outer policies
\[
\tpsi^{\kk}(z) = \frac{1 + z}{2} 
\exp\!\left(\gamma_\kappa\, \frac{z-1}{1 +z+k^{-1}}\right).
\]
Then $\{\tpsi^{\kk}\colon k \in \mathbb{N}\}$ is an $\epsilon$-optimal sequence, that is, $\displaystyle {\cal C}^* = \lim_{k \to \infty} {\cal C}(\tpsi^{\kk})$.

\end{enumerate}
\end{thm}

For $\kappa \ge \sqrt{5}$, the optimal policy takes a simple and implementable form: the retailer’s orders become a convex combination of current and one-period lagged demand,
\[
O_t = \tpsi^*_0\,D_t + (1 - \tpsi^*_0)\,D_{t-1}.
\]
In contrast, for $\kappa < \sqrt{5}$, the $\epsilon$-optimal sequence $\{\tpsi^{\kk}\}$ has no finite-order ARMA representation, and the limiting transfer function
\begin{equation}\label{eq:limit_policy}
\tpsi^\infty(z) := \lim_{k \to \infty}\tpsi^{\kk}(z)
= \frac{1 + z}{2}\,\exp\!\left(\gamma_\kappa\,\frac{z-1}{1 + z}\right)
\end{equation}
contains an exponential singular inner factor. Let
\[
\sigma^2_{\I^*} := \lim_{k \to \infty} \sigma^2_{\I}(\tpsi^{\kk})
\qquad \mbox{and}\qquad 
\sigma^2_{\Su^*} := \lim_{k \to \infty} \sigma^2_{\Su}(\tpsi^{\kk}).
\]
For $\kappa < \sqrt{5}$, the retailer’s inventory variance is continuous at the limit, $\sigma^2_{\I^*} = \sigma^2_{\I}(\tpsi^\infty)$, but the supplier’s MSFE exhibits an upward jump, with $\sigma^2_{\Su^*} < \sigma^2_{\Su}(\tpsi^\infty)$. Hence ${\cal C}^* < {\cal C}(\tpsi^\infty)$ and the limiting policy $\tpsi^\infty$ is strictly suboptimal.

The presence of the singular inner factor in \eqref{eq:limit_policy} reflects the intrinsic complexity of the optimization problem when $\kappa$ is small (i.e., when the supplier costs are relatively more important). In particular, the $\epsilon$-optimal policies $\{\tpsi^{\kk}\}$ do not admit a finite-order ARMA representation, which makes them impractical to implement in practice.

Motivated by these challenges, the remainder of the paper proposes a new perspective for interpreting the retailer’s inventory policy and develops a practically implementable approximation framework; together, these developments constitute our main contributions:
\begin{enumerate}
\item Conceptually, our fundamental contribution is our interpretation of an optimal policy through the lens of \emph{information delay}: the supplier’s limited access to real-time demand information forces the retailer to strategically postpone replenishment, using delay as a mechanism to smooth and increase the predictability of its orders, thereby trading off its own inventory costs against the supplier’s forecast accuracy. Importantly, the optimal policy is not characterized by \emph{pure} delay alone, but requires the inclusion of \emph{group delay} (as defined later), which captures how the timing of information is redistributed across lags. This perspective clarifies how information delay and group delay jointly shape the structure of the optimal policy and the distribution of costs across the supply chain.

\item We construct a concrete, implementable ARMA approximation to the limiting policy $\tpsi^\infty(z)$ that attains near-optimal cost while remaining rational, computationally tractable, and transparent. This yields a practically deployable policy for decentralized supply chains.

\item We extend the modeling framework to account for situations in which the supplier’s forecast is based not on the entire order history, but on a finite window of past observations. This modification reflects more realistic information structures in practice and introduces a new dimension to the analysis. We establish convergence results showing that as the length of the supplier’s memory increases, the associated cost approaches the infinite-memory benchmark, thereby bridging the gap between the idealized full-memory setting and implementable, finite-memory  forecasting policies. 
\end{enumerate}
An important by-product of developing our information-delay perspective is the introduction of a new set of mathematical tools and models to study the connection between inventory management policies and their forecastability in the context of supply chains.

\section{Information Delay}\label{sec:Delay}

In this section we provide an interpretable characterization of the retailer’s inventory policy through the lens of \emph{information delay}, formalized via the  \emph{group delay} of the policy's transfer function (Definition~\ref{Group Delay Def}). We show that the delay embedded in a policy acts as a signaling channel that the supplier can exploit to improve order forecasts, but that it comes at the cost of higher inventory volatility and delayed replenishment for the retailer.

\begin{figure}[hbt]
   \begin{center}
\includegraphics[width=19.7cm]{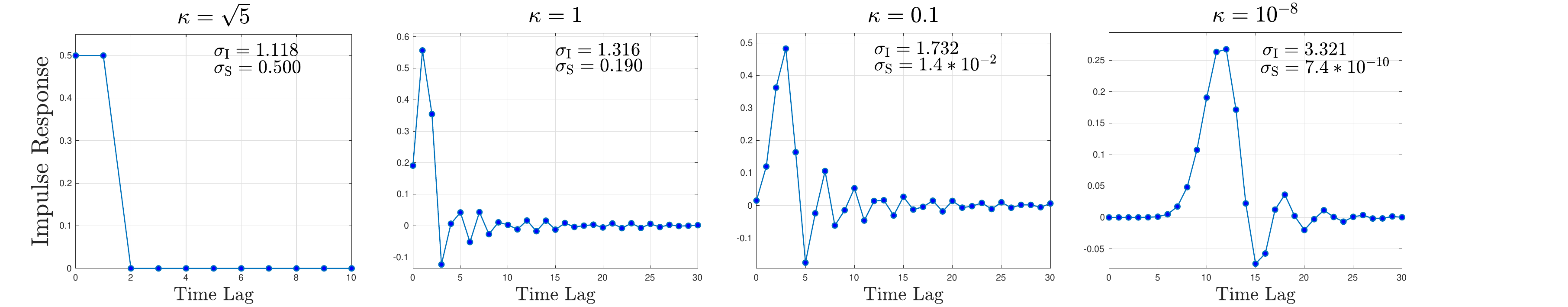}
   \end{center}
   \caption[Model Description]{\footnotesize Impulse–response coefficients $\{\hpsim_{n}\}_{n\ge 0}$ for different values of  $\kappa$ and $m=100$.}
   \label{Fig:ImpulseCoef}
\end{figure}

Before turning to details, we illustrate the notion of information delay in Figure~\ref{Fig:ImpulseCoef}. The figure shows the impulse–response coefficients $\{\hpsim_{n}\}_{n\ge 0}$ of our $\epsilon$-optimal ARMA policy $\hpsim(z)=\sum_{n=0}^\infty \hpsim_{n} z^n$ (detailed in section \ref{sec: Optimal Orders}) for several values of $\kappa$, together with the corresponding values of $\sigma_{\I}(\hpsim)$ and $\sigma_{\Su}(\hpsim)$. In the time domain, under the bounded–inventory constraint $\hpsim(1)=\sum_{n=0}^\infty \hpsim_n=1$, we interpret $\hpsim_{n}$ as the fraction of market demand $D_t$ realized at time $t$ that appears in the order $O_{t+n}$ placed at time $t+n$. Thus, the sequence $\{\hpsim_{n}\}_{n\ge 0}$ characterizes the retailer’s delayed replenishment response to demand %\prem{Maybe change figure with higher value of $m$}.
%\cliff{The discussion in this paragraph, with slight modification, would apply to any policy $\tpsi$. So maybe this discussion should be moved earlier and presented for general $\tpsi$.}

At one extreme, for $\kappa=\sqrt{5}$ (left panel), the impulse response exhibits minimal delay, distributing the demand impact evenly between lags 0 and 1. At the other extreme, for $\kappa=10^{-8}$ (right panel), the response is substantially delayed, with a noticeable reaction to demand only emerging around lag 6. As the delay increases (left to right), the supplier’s root MSFE $\sigma_{\Su}$ decreases while the retailer’s inventory volatility $\sigma_{\I}$ rises, illustrating how optimal delay governs the trade-off between retailer and supplier inventory costs. We quantify this relationship in Theorems \ref{thm:GD_MSFE_Improved} and  \ref{Inventory_Breakdown}, which connect $\sigma_{\I}(\tpsi)$ and $\sigma_{\Su}(\tpsi)$ to the information delay embedded in $\tpsi$.

We adopt the signal-processing concept of \emph{group delay} to formalize information delay, specializing it to the class of admissible inventory policies $\tpsi \in \widetilde{\Psi}$, analytic functions $\tpsi(z)=\sum_{k=0}^{\infty} \tpsi_k z^{k}\in\mathbb{H}^2$ with real MA coefficients ${\tpsi_k}$ that satisfy the bounded-inventory normalization $\tpsi(1)=1$.

\begin{definition}
\label{Group Delay Def}{\rm \bf (Group Delay)}
Let $\tpsi(z) = \sum_{k = 0}^{\infty} \tpsi_k z^{k} \in \mathbb{H}^2$ be a $z$-transform with Taylor coefficients \mbox{$\{\tpsi_k\} \subseteq \mathbb{R}$} and $\tpsi(1)=1$. We define the \emph{group delay} of $\tpsi$ as
\begin{equation*}
\text{Group Delay: } \GD(\tpsi) \equiv \tpsi'(z)\Big|_{z=1} = \sum_{k=0}^{\infty} k \,\tpsi_k.
\end{equation*}
\end{definition}
In words, the quantity $\sum_{k=0}^{\infty} k\,\tpsi_k$ gives a time-domain view of group delay by measuring the “center of mass” of the impulse response (since $\tpsi(1)=1$ implies $\sum_{k=0}^{\infty} \tpsi_k=1$). A larger group delay shifts, on average, the future time at which a current demand shock is fully reflected in the output signal (i.e., retailer's orders). To further build intuition for this definition in our supply-chain setting, let us review a few concrete  examples.
\vspace{0.1cm}

{%\setstretch{1.1}
\begin{exm}\label{exm-delaytype}~~{\sf 
\vspace{0.2cm}

\begin{enumerate}[{\rm (a)}]

\item {\sc Pure Delay}. Under a  pure delay the retailer's orders satisfy $O_t = D_{t-k}$ for some integer $k \in \mathbb{N}$, that is, the retailer places orders by simply delaying its demand by exactly $k$ periods. In this case, the $z$-transform of the retailer's orders is $\tpsi(z) = z^k$ and its group delay is given by $\GD(\tpsi) = k$, which is precisely what we expect a measure of delay should be in this setting. \vspace{0.2cm}

\item {\sc Non-Pure Delay}. Suppose the retailer replenishes its inventory according to the MA(1) orders \mbox{$O_t = \tpsi_0 D_t + \tpsi_1 D_{t-1}$}, where $\tpsi_0 + \tpsi_1 = 1$. In this case, $\tpsi(z)= \tpsi_0+ (1-\tpsi_0)\,z$ and the group delay of these orders equals $\GD(\tpsi) = 1-\tpsi_0$. For instance, if $\tpsi_0  = 0.5$, then $\GD(\tpsi) = 0.5$, which can be interpreted as the retailer partially delaying the replenishment of its demand by half a period on average. 
It is worth noting that if $\tpsi_0 > 1$, then $\GD(\tpsi)  < 0$, indicating that the retailer's orders exhibit a negative delay in this case. Intuitively, $\tpsi_0 > 1$ implies that the contribution of the retailer's demand $D_t$ in period $t$ to the order $O_t$ placed in the same period exceeds 100\%, and the excess must therefore be ``returned'' in future periods. Thus,  group delay, as defined in \cref{Group Delay Def}, captures the effective timing of how demand in a given period impacts the retailer's orders in subsequent periods. %We will see in the discussion below, why this strategy is suboptimal in our supply chain. 
\vspace{0.2cm}

\item {\sc Inner Functions:} An important special class of transfer functions for which group delay is guaranteed to be nonnegative is the class of inner functions (all-pass filters) in $\mathbb{H}^2$. Let $I(z) \in \mathbb{H}^2$ be an inner function, and let the polar decomposition of $I$ on the unit circle $\mathbb{T}$ be given by $I(e^{-i\lambda}) = r(\lambda)e^{-i\phi(\lambda)}$, where $r(\lambda) > 0$ denotes the gain and $\phi(\lambda)$ the phase. Then, the group delay of $I$ satisfies $\GD(I) = \phi'(0)$. That is, in a frequency-domain sense, the group delay of an inner function quantifies the effective time lag experienced by signal components near zero frequency as they pass through the system, capturing how the phase of the output signal varies with frequency at the origin.
The fact that $\GD(I) = \phi'(0) \geq 0$ follows from the  Julia-Wolff-Carathéodory theorem, with strict inequality for whenever $I$ is nonconstant; see Theorem 2.4 in \citep{dang2011analytic} ~$\diamond$
\end{enumerate}}
\end{exm}}

In the context of our supply chain system, one might intuitively expect that if the retailer’s orders exhibit positive group delay, then this delay should enhance the supplier’s ability to forecast those orders. However, this intuition is generally false; the class of inner functions provides a canonical counterexample. As noted in Example~\ref{exm-delaytype}(c), inner functions have nonnegative group delay but do not alter the spectral density of the order process. Since the supplier’s forecastability, as captured by the MSFE in Kolmogorov’s formula~\eqref{eq:Kolmogorov}, depends only on this spectral density, an all-pass component leaves the supplier’s forecast error unchanged while increasing the retailer’s inventory cost. The following simple example makes this point concrete.\vspace{0.2cm}
\vspace{0.2cm}

{%\setstretch{1.1}
\begin{exm}[Pure delay continued]\label{delaytype-inventory}{\sf 
Consider the inner function $\tpsi(z)=I_k(z)=z^k$ with $k\in\mathbb{N}$ in Example~\ref{exm-delaytype}(a); it has group delay $\GD(I_k)=k$. By~\eqref{eq:Kolmogorov}, the supplier's MSFE satisfies $\sigma^2_{\Su }(I_k)=1$, independent of $k$. By contrast, by~\eqref{eq_SigmaI_z}, the retailer’s inventory variance grows linearly, $\sigma^2_{\I}(I_k)=1+k$. Intuitively, delaying orders by $k$ periods has no effect on the supplier's forecast error but exposes the retailer to the cumulative uncertainty of $k$ independent demand shocks, raising the retailer's cost ${\cal C}^{\R}$ by a factor $\sqrt{1+k}$ (recall that the demand-shock variance is normalized to one). Thus, a pure delay is strictly harmful to the retailer and the supply chain as a whole. ~$\diamond$}
\end{exm}
}

The following result establishes a direct link between the group delay of a policy $\tpsi(z)$ and the supplier’s mean squared forecast error (MSFE), and demonstrates that an invertible order \textit{must} possess positive group delay in order to keep the supplier's forecasting error low.

\begin{thm}\label{thm:GD_MSFE_Improved}
Consider a policy $\tpsi \in \widetilde{\Psi}$ with factorization $\tpsi(z) = Q(z)\,I(z)$,
where $Q(z)$ is an outer function, and $I(z)$ is an inner function. Suppose:
\begin{align*}
    \int_{-\pi}^{\pi} \frac{\log |Q(e^{-i\lambda})|}{|e^{i\lambda} - 1|^2} d\lambda < \infty
\end{align*}
Then:
\begin{equation}
\label{eq: Forecasting Two-Sided}
c e^{-\,\GD(Q)} \;\geq\; \sigma^2_{\Su}(\tpsi) \;\geq\; c^{-1}e^{-\,\GD(Q)} \geq c^{-1}e^{-\,\GD(\tpsi)}
\end{equation}
where:
\begin{equation}
\label{eq: Constant Definition}
    c \equiv \frac{1}{2\pi}\int_{-\pi}^{\pi} \frac{\cos \lambda}{1 - \cos \lambda} \log |Q(e^{i\lambda})|d\lambda  < \infty
\end{equation}
\end{thm}

This theorem demonstrates that for a policy $\tpsi$ to achieve a small MSFE, it is both necessary and sufficient (under weak regularity conditions) for its outer factor $RQ$ to exhibit a \emph{large} group delay. In other words, $\GD(Q)$, and consequently $\GD(\tpsi)$, sets an \emph{exponential upper/lower bound} on the supplier's attainable MSFE. In our applications, the value $c$ in \eqref{eq: Constant Definition} can typically be considered a constant fixed over a family of $\tpsi$; for example, when $Q(z) = \prod_{k = 1}^n  \frac{z - e^{i\theta_k}}{1 - e^{i\theta_k}} \prod_{l = 1}^m \Big(\frac{z - a_l}{1 - a_l}\Big)^{\alpha_l}$ (where $\alpha_1 \in \mathbb{Z}$) and all $|a_l| > 1$ and sufficiently separated from $1$.

\begin{cor}
\label{Constant is Controlled}
Let $Q(z) = \prod_{k = 1}^n  \frac{z - e^{i\theta_k}}{1 - e^{i\theta_k}} \prod_{l = 1}^m \Big(\frac{z - a_l}{1 - a_l}\Big)^{\alpha_l}$ where $\alpha_1 \in \mathbb{Z}$ and $a_l, e^{i\theta_k} \in \bar{\mathbb{D}}^c \cap \mathbb{B}(1, \delta)^{c}$ for some $\delta > 0$. Then:
\begin{equation*}
    |c| \leq C(\delta)\Big(n + \sum_{l = 1}^{m} |\alpha_l|\Big)
\end{equation*}
where $c$ is as defined above and $C(\delta)$ is a constant depending only on $\delta$.
\end{cor}

While Theorem \ref{thm:GD_MSFE_Improved} suggests that group delay is a necessary requirement for lowering the supplier’s costs, we expect group delay to increase costs for the retailer. Indeed, as noted in Example \ref{delaytype-inventory}, group delay exposes the retailer to uncertainty, forcing him to accumulate inventory in order to hedge against instantaneous demand shocks. Our next result formalizes this notion, and provides a complete, explicit breakdown of the retailer's costs with respect to the invertible (outer) and all-pass (inner) components of the order. Indeed, Theorem \ref{Inventory_Breakdown} demonstrates that applying an all-pass filter to an invertible order (i.e. inserting an inner factor into its transfer function) invariably increases the retailer’s inventory variance $\sigma^2_{\I}$ by \textit{exactly} its group delay.

\vspace{0.2cm}

\begin{thm}\label{Inventory_Breakdown}
Let $\tilde{\psi}(z) = Q(z)B(z)S(z)$ be expressed in terms of its Beurling factorization, where 
\[ B(z)  = \prod_{n = 1}^{\infty} \frac{a_n}{\bar{a}_n} \frac{a_n - z}{1 - \bar{a}_n z} \qquad \mbox{and}\qquad
    S(z)  = \exp\Big(-\int \frac{\xi - z}{\xi + z}d\mu(\xi)\Big).
\]
Suppose:
\begin{equation}
\label{eq: Integrability Condition}
    \sum_{n = 1}^{\infty} \frac{1 - |a_n|^2}{|1 - a_n|^2} + \int_{-\pi}^{\pi} \frac{d\mu(\xi)}{|\xi - 1|^2} + \int_{-\pi}^{\pi} \frac{\log |Q(e^{-i\lambda})|}{|e^{i\lambda} - 1|^2} d\lambda < \infty
\end{equation}
with $|a_n| < 1$ for all $n \in \mathbb{N}$ and $\mu \perp d\lambda$ \footnote{The symbol $\perp$ indicates two measures are mutually singular.}. Then, 
\begin{equation}
\label{eq: Inventory Breakdown Equation}
\sigma^2_{\I}(\tpsi) = 1 + \GD(Q) + \GD(B) + \GD(S) - \frac{1}{2}\partial_r f_{Q}(1)
\end{equation}
where  $\displaystyle f_{Q}(re^{-i\lambda}) = b_0 + 2\sum_{k \geq 1} b_k r^k \cos k\lambda$ with $\{b_i\}_{i \geq 0}$ the autocovariances of $\tilde{\psi}$ and $\partial_r$ denotes the radial derivative. Furthermore, 
\begin{equation*}
    \sigma^2_{\I}(\tpsi) = \sigma^2_{\I}(Q) + \GD(B) + \GD(S) \geq \sigma^2_{\I}(Q).
\end{equation*}
\end{thm}
Thus, inner factors $I(z)=B(z)S(z)$ increase the retailer’s inventory volatility by exactly their total group delay, so that $\sigma^2_{\I}(\tilde{\psi})-\sigma^2_{\I}(Q)=\GD(I)=\GD(B)+\GD(S)$. Moreover, as in Example~\ref{exm-delaytype}(a), the group delay of $I$ equals the derivative at zero of its phase $\phi(\lambda)=\arg I(e^{-i\lambda})$, namely $\GD(I)=\phi'(0)$. This derivative is the slope of the best linear approximation of the phase near frequency zero and can be interpreted as approximating $I(z)$ by the ``pure'' delay $z^{\phi'(0)}$.  Hence $\GD(I)=\phi'(0)$ serves as a proxy for the effective time delay induced by the all-pass (inner) filter $I$.

Furthermore, from \eqref{eq: Inventory Breakdown Equation}, we see that outer factor $Q$ influences $\sigma^2_{\I}$ through \textit{both} its group delay $\GD(Q)$ and the quantity $\partial_r f_Q(1)$, which depends only on the order spectral density and measures the ``center'' of mass of its autocorrelation profile (indeed, $\partial_r f(1) = \sum_{k \geq 1} kb_k$ by definition). Intuitively, we can see $|\partial_r f_Q(1)|$ as measuring the failure of the factor $Q$ to be inner (indeed if $|Q(e^{i\theta})| = 1$ a.e. on $[-\pi, \pi)$, then $\partial_r f_Q(1) = 0$). This observation will be critical in Proposition \ref{Singular Inner Approximation}, where we will see $|\partial_r f_{Q_m}(1)| \to 0$ for a sequence $\{Q_m\}_{m \geq 1}$ approximating the all-pass component of $\tpsi^\infty(z)$. 

We note that Theorem~\ref{Inventory_Breakdown} significantly strengthens the results in \cite{caldentey2024information}, in particular Lemma 3 therein, by showing that \emph{any} inner factor (not only a Blaschke product) is suboptimal, as it raises the retailer’s inventory volatility without reducing the supplier’s MSFE. Because inner factors are spectrally invisible to the supplier (they leave the retailer’s order spectrum unchanged), they yield no forecasting benefit. Thus, to improve forecastability, the retailer must encode delay through an \emph{invertible} (outer) order. Finally, we remark that the weak regularity assumptions in \eqref{eq: Integrability Condition} are made to ensure the finiteness of $\sigma^2_{\I}$ and are always satisfied for rational $Q$ with real roots/poles, atomic $\mu$, and finite Blaschke products $B$ (the setting considered in this paper). 

%\prem{Do we need to discuss the properties of ${\tpsi_k}$? It's not a focus of the paper --%- I would get rid of this paragraph}. Interestingly, and somewhat counterintuitively, the %sequence of $\eps$-optimal policies ${\tpsi_k}$ from Theorem~\ref{thm_CGHZ} converge to %the limit policy $\tpsi^\infty(z)$ in \eqref{eq:limit_policy}, whose transfer function %contains the singular inner factor $I_\infty(z)=\exp\big(\gamma_\kappa\,(z-1)/(1+z)\big)$. %This suggests that the sequence ${\tpsi_k}$ is effectively approximating the informational-%delay benefits of an all-pass inner filter using outer (invertible) ones 

Taken together, these observations motivate using pure-delay filters as a \emph{blueprint} for designing outer (invertible) filters that approximate the same informational-delay benefits. Indeed, the combination of Theorems \ref{thm:GD_MSFE_Improved} and \ref{Inventory_Breakdown} demonstrates that a \textit{positive} group delay in $\tpsi$ is needed to optimize total costs --- while positive group delay increases costs for the retailer \textit{linearly} (Theorem \ref{Inventory_Breakdown}), under weak regularity conditions, decreases costs for the supplier \textit{exponentially} (Theorem \ref{thm:GD_MSFE_Improved}). Moreover, Proposition ~\ref{Singular Inner Approximation} below shows that a broad class of (all-pass) singular inner filters $S$ (which impose fractional delay as in $\tpsi^{\infty}$ in \eqref{eq:limit_policy}) can be realized as the limit of invertible filters whose phase derivatives remain strictly positive, thereby preserving the informational role of group delay while making the delay \emph{detectable} to the supplier. In designing our approximation scheme, we look to Theorem \ref{Inventory_Breakdown} --- we choose our sequence $\{\tilde{Q}_m\}_{m \geq 1}$ of invertible filters so that $\text{GD}(\tilde{Q}_m) = \text{GD}(S)$ $\forall m \geq 1$ and with $\partial_r f_{\tilde{Q}_m}(1) \to 0$ as $m \to \infty$.

\begin{proposition}
\label{Singular Inner Approximation}
Let $Q(z)\in\mathbb{H}^2$ and fix an atomic singular inner function $S(z)\in\mathbb{H}^2$ \footnote{A singular inner function whose singular measure $\mu$ is atomic, see statement of Theorem \ref{Inventory_Breakdown} for definition of $S$.} with boundary values $S(e^{-i\lambda})=e^{-i\phi(\lambda)}$ a.e.\ on $\mathbb{T}$, that satisfies the bounded inventory condition $S(1)=1$ (i.e., $\phi(0)=0$) and the integrability condition\footnote{The condition \eqref{eq: Decay Condition on Phase} characterizes the rate at which $\phi(\lambda) \to 0$ as $\lambda \to 0$, requiring $\phi^2(\lambda) \asymp \tan \frac{\lambda}{2}$; this is comparable to requiring $\sigma^2_{\I}(S) < \infty$ (note that $|1 - z|^2 = \sin^2 \frac{\lambda}{2}$ in \eqref{eq_SigmaI_z}), where $\asymp$ denotes asymptotic equivalence up to positive multiplicative constants.
}. Let $g_Q(e^{-i\lambda}) = |Q(e^{-i\lambda})|^2$ and suppose that:
\begin{equation}
\label{eq: Decay Condition on Phase}
    \int_{-\pi}^{\pi}  |g'_Q(e^{i\lambda})| \phi^2(\lambda) \cot \frac{\lambda}{2} d\lambda + \int_{-\pi}^{\pi} g_Q(e^{i\lambda}) \Big|\phi(\lambda)\phi'(\lambda)\cot \frac{\lambda}{2}\Big| d\lambda < \infty
\end{equation}
For some some $m \in \mathbb{N}$, consider the approximation $\tilde{Q}_m$ to $S$ given by:
\begin{equation*}
    %\tilde{Q}(e^{-i\lambda}) = \Big(1 - \frac{i\log S(e^{-i\lambda})}{m}\Big)^m
    \tilde{Q}_m(z)=\left(1-{\log(S(z))\over m}\right)^{-m},
\end{equation*}
Then, $\tilde{Q}_m$ is outer and satisfies $ \cal{GD}(\tilde{Q}_m) = \cal{GD}(S) = \phi'(0)$. Moreover,  we have that 
\begin{equation}
\label{eq: Asymptotic Decomposition}
    \sigma^2_{\I}(Q\tilde{Q}_m) = \sigma^2_{\I}(Q) + \cal{GD}(\tilde{Q}_m) + \mathcal{O}\Big(\frac{1}{m}\Big)=\sigma^2_{\I}(QS)+\mathcal{O}\Big(\frac{1}{m}\Big).
\end{equation}

\end{proposition}
Hence, the outer function $\tilde{Q}_m$ approximates the delay mechanism of $S$ by approximating the phase $\phi(\lambda)$ of $S$. Indeed, from the proof of the proposition 
one can see that the phase $\tilde{\phi}_m$ of $\tilde{Q}_m$ satisfies $\tilde{\phi}_m(0)=\phi(0)$ and $\tilde{\phi}_m'(0)=\phi'(0)$.\vspace{0.2cm} 

In section \ref{sec: Optimal Orders}, we will apply Proposition \ref{Singular Inner Approximation} to build an invertible ARMA class of orders that approximate the optimal policy. Our approach hinges on replacing the singular inner factor $S(z)=\exp\bigl(\gamma_\kappa\mathbb{}(z-1)/(1+z)\bigr)$ of the limiting policy $\tpsi^\infty(z)$ with the outer factor $\tilde{Q}_m(z)=\bigl(1-\tfrac{1}{m}\log S(z)\bigr)^{-m}$.

\section{ARMA Approximations}
\label{sec: Optimal Orders}
%\yichen{I actually have another question to this ARMA, may be not highly relevant, we can discuss separately. If the order (the demand of the supplier) is ARMA and Gaussian, the process than be represented using kalman filter with state space degree $m$, and thus it is a finite-memory Hidden Markov Model. I am not sure if that will force us to compare our result with those earlier development based on markov properties. I guess I need to ask what are the state-of-art results on this direction?} \cliff{My understanding is that in this setting, given that the ARMA parameters are known, and given that the supplier's forecast uses a finite past, the minimum-MSFE forecast is the solution to the Yule-Walker equations, which is what we are using. If you are asking what the optimal forecast is based on a finite past and unknown ARMA parameters, I believe this is an open probem. Maybe you could try to attack it using Hidden Markov Models, but I would be surprised if you could obtain an optimal solution using HMMs.}
 As discussed in section \ref{sec:preliminaries},  for small values of $\kappa$,  the limiting policy  involves a non-rational transfer function with a singular inner factor. In this section, we leverage Proposition \ref{Singular Inner Approximation} to introduce a class of $\epsilon$-optimal invertible ARMA orders that we use to approximate this solution. We numerically evaluate the performance of this class and show that even low-order ARMA models achieve near-optimal performance, striking a practical balance between accuracy and tractability while providing a compact and interpretable representation of the optimal information-delay mechanism discussed in the previous section. We also derive asymptotic bounds for the suboptimality that quantify the rate of convergence of this class of policies as their ARMA degree grows large. Finally, from Theorem~\ref{thm_CGHZ}, we already know that the optimal solution to the optimization problem in \eqref{eq:Cost_SC_z} is a simple MA(1) policy when $\kappa \geq \sqrt{5}$. Thus, we concentrate the analysis on the case $\kappa < \sqrt{5}$.

Our proposed approximating procedure takes the limiting policy 
$$\tpsi^\infty(z) =\frac{1 + z}{2} \exp\!\left(\gamma_\kappa\, \frac{z-1}{1 + z}\right),$$ as a target reference and replaces the singular inner factor 
using  Euler's limit representation of \mbox{$\exp(-x) \approx (1+x/m)^{-m}$} for large $m \in \mathbb{N}$.  Using this approximation scheme, we obtain the following sequence of causal, \textit{invertible} ARMA orders $\hpsi^{m}$ to approximate $\tilde{\psi}^\infty$ 
\begin{align}\label{eq:tpsi_m}
        \hpsi^{m}(z) & \equiv \frac{1 + z}{2}\Bigg(1 + \frac{\gamma_\kappa(1 - z)}{m(1 + z)}\Bigg)^{-m}  = \underbrace{\Big(\frac{1 + z}{2}\Big)^{m + 1}}_{\text{MA}} \cdot \underbrace{\Big(\frac{1 + \frac{\gamma_\kappa}{m}}{2}\Big)^{-m}\Big(1+\frac{1 - \frac{\gamma_\kappa}{m}}{1 +\frac{\gamma_\kappa}{m}}z\Big)^{-m}}_{\text{AR}},
\end{align}
where, in the second equality, we decompose $\hpsi^{m}$ into an MA component of order $m{+}1$ and an AR component of order $m$, each satisfying the bounded–inventory condition. Thus, from a practical standpoint, the parameter $m$ captures the ``{\em complexity}'' of our proposed approximating policy $\hpsi^{m}(z)$.

As shown below, for fixed $m\in\mathbb{N}$, the MA part helps the supplier by damping forecasting error, while the AR part helps the retailer by reducing group delay (see Section~\ref{sec:Delay}).  Increasing $m$ improves the retailer’s approximation to $\tilde{\psi}^\infty$ but harms the supplier by amplifying the MSFE $\sigma^2_{\Su}$ (see the proof of Theorem \ref{thm:GD_MSFE_Improved}).

In passing, we note that although other rational-approximation schemes for the singular inner factor in $\tpsi^\infty(z)$ are possible, our proposed scheme has the key property of matching the group delay of the all-pass component of $\tpsi^\infty(z)$, which makes it a natural choice for both analysis and implementation. This follows directly from Proposition~\ref{Singular Inner Approximation}.

\begin{cor}\label{cor:ARMAapprox} The approximating ARMA policies $\hpsi^m(z)$ in \eqref{eq:tpsi_m} preserve the group delay of the limiting policy $\tpsi^\infty(z)$, i.e., $\displaystyle \GD(\hpsi^m)=\GD(\tpsi^\infty)={\D \over \D z} \tpsi^\infty(z)\Big|_{z=1}={1 +\gamma_\kappa \over 2}$.
\end{cor}

Recall from Theorem \ref{thm_CGHZ} that $\gamma_{\kappa}$ is the unique solution to $\kappa^2=(5+2\,\gamma)\,\exp(-2\,\gamma)$; thus as $\kappa \to 0$, $\gamma_{\kappa} \to \infty$.  Hence, Corollary \ref{cor:ARMAapprox} illustrates that prioritizing the supplier ($\kappa := K^{\R} / K^{\Su} \ll 1$) mandates increasing the delay $\GD(\hpsi^m)$ of our order policy, which further confirms the conclusion of Theorem \ref{thm:GD_MSFE_Improved} that large delay is \textit{necessary} to achieve low inventory costs for the supplier.

To assess the performance of the approximating policies in \eqref{eq:tpsi_m}, \cref{Table_approx} reports the relative cost ${\cal C}(\hpsi^{m})/{\cal C}^*$ as a function of $m$ and the retailer’s relative cost weight $\kappa$. Overall, the approximation performs uniformly well across values of $\kappa$ for $m \ge 10$. Moreover, for $\kappa \ge 1$, that is, when the retailer’s cost carries greater weight than the supplier’s, the approximation already delivers excellent performance for $m \ge 1$, indicating that even a modest amount of order complexity can yield substantial cost reductions for the supply chain.

\setlength{\tabcolsep}{10pt} % Increase column spacing globally
\renewcommand{\arraystretch}{1.1} % Optional: adds vertical padding
\begin{table}[h!]
\centering
\rowcolors{4}{gray!10}{white} % start from data row 1
\begin{tabular}{c| c c c c c c c c}
\multicolumn{1}{c}{} & \multicolumn{8}{c}{$m$} \\
\cmidrule(lr){2-9}
$\kappa$ &0 & 1 & 2 & 5 & 10 & 20 & 50 & 100 \\
\midrule
0.001  & 204.707 & 10.767 & 2.824 & 1.271 & 1.111 & 1.046 & 1.014 & 1.006 \\
0.01 & 23.390   & 2.255  & 1.446 & 1.139 & 1.059 & 1.025 & 1.008 & 1.003 \\
0.1 & 3.249   & 1.284  & 1.135 & 1.047 & 1.020 & 1.008 & 1.003 & 1.001 \\
0.5 &  1.298  & 1.052  & 1.027 & 1.009 & 1.004 & 1.002 & 1.001 & 1.000 \\
1.0  & 1.074  & 1.012  & 1.006 & 1.002 & 1.001 & 1.000 & 1.000 & 1.000 \\
2.0 &  1.001   & 1.000  & 1.000 & 1.000 & 1.000 & 1.000 & 1.000 & 1.000 \\
$\sqrt{5}$ & 1.000 & 1.000  & 1.000 & 1.000 & 1.000 & 1.000 & 1.000 & 1.000 \\
\bottomrule
\end{tabular}\vspace{0.1cm}

\caption{Relative cost performance ${\cal C}(\hpsi^{m})/{\cal C}^*$ of the approximating policy $\hpsi^{m}$ as a function of $\kappa$ and $m$.}\label{Table_approx}
\end{table}
\vspace{0.3cm}

{%\setstretch{1.1}
\begin{exm}{\sf 
To further illustrate the results in Table~\ref{Table_approx}, 
consider the case $\kappa = 1$, where the retailer’s and supplier’s costs 
contribute equally to the total supply chain cost. 
In this case, choosing $m = 1$ yields a simple ARMA(1,2) ordering policy:
\[
O_t = \frac{\gamma_\kappa - 1}{1 + \gamma_\kappa}\,O_{t-1} 
+ \frac{1}{2(1 + \gamma_\kappa)}\,(D_t + 2D_{t-1} + D_{t-2}),
\qquad \text{with } \gamma_\kappa \approx 0.968,
\]
whose cost is only $1.2\%$ higher than that of an optimal policy. Increasing the order to $m = 2$ (i.e., an ARMA(2,3) policy) further reduces this gap in half.}~$\diamond$
\end{exm}}\vspace{0.2cm}

The next result reinforces the results in Table~\ref{Table_approx} by providing an analytic quantification of the asymptotic performance of $\hpsi^m(z)$ as $m$ grows large. 

\begin{proposition}
\label{Inventory Approximation} For any $m \in \mathbb{N}$ and for some constant $C_\kappa > 0$ depending only on $\gamma_\kappa$,
\begin{enumerate}[label={\rm (\roman*)}, leftmargin=3em]
\item $|\sigma^2_{\I}(\tpsi^\infty) - \sigma^2_{\I}(\hpsi^{m})| \le C_\kappa\, m^{-1}$, 
\label{Inventory Error}
\item  $|\sigma^2_{\Su*} - \sigma^2_{\Su}(\hpsi^{m})| \leq C_\kappa\, m^{-1}$, and\label{Supplier's Error}
\item $|{\cal C}(\hpsi^{m})-{\cal C}^*| \leq C_\kappa\, m^{-1}$. \label{Total Cost Error}
\end{enumerate}
\end{proposition}
\cref{Inventory Approximation} highlights the trade-off between optimality and implementability inherent to our proposed class of approximating policies $\{\hpsi^m\}$ in \eqref{eq:tpsi_m}. On one end, from a practical standpoint, we would like to choose a small value of $m$ to reduce the complexity of the policy. On the other hand, we would like to make $m$ as large as possible to reduce the policy's optimality gap.

\subsection{Finite Memory}\label{sec:finitepast}
In this section, we extend the results in the previous section by 
investigating the inventory manager’s problem when the supplier has finite memory and builds its forecasts on a finite window of past orders rather than the entire order history. In practice, the supplier’s reliance on a limited order history may arise from suboptimal forecasting methods or computational constraints, leading to a finite-memory formulation that captures realistic informational limits and highlights a natural trade-off between tractability and performance.
 Our analysis focuses on bounding the performance gap between such finite-memory policies and the idealized infinite-history benchmark, thereby quantifying the cost of limited information. In addition, we characterize the rate at which this gap vanishes as the supplier’s memory length increases, showing that the costs of finite-memory policies converge to those of the full-history setting at an explicit rate. This allows us to assess how quickly implementable policies approach the theoretical benchmark as the horizon of past observations expands. In addition, in contrast to \cref{Inventory Approximation}, we show that when the supplier has finite memory, increasing the complexity of the retailer's policy might be counterproductive from an optimality standpoint.

We consider a situation in which at any moment in time the supplier uses only the last $n+1$ orders, $\{O_t, O_{t-1}, \ldots, O_{t-n}\}$, to forecast the retailer’s next-period order. In this case, $\mathcal{F}^{\Su}(t) = \mathcal{F}_t^{\ti{O},n}$, where $\mathcal{F}_t^{\ti{O},n}$ is the $\sigma$-algebra generated by this finite history of the retailer’s orders.  With a finite memory, the supplier constructs the best linear $n$-past predictor $\widehat{O}^n_{t+1}$ of the retailer's next period orders  by computing the projection of $O_{t+1}$ into the linear span generated by $\{1,O_t, O_{t-1}, \ldots, O_{t-n}\}$. That is, 
    $$\widehat{O}^n_{t+1}=c+\sum_{k=0}^n c_k\,O_{t-k},$$
where the coefficients $\mathfrak{c}=(c,c_{0}, \ldots, c_{n})$ realize the minimum in the following optimization problem:
\begin{equation}
\label{n-past MSFE definition}
\sigma^2_{\Su,n} \equiv \min_{\mathfrak{c}} \mathbb{E}\Big[\Big(O_{t+1} - c - \sum_{k = 0}^n c_{k}\, O_{t-k}\Big)^2\Big].
\end{equation}
Note that, by stationarity, the quantity $\sigma^2_{\Su,n}$, the \textit{$n$-past prediction error}, is independent of $t$. Also, in general, we have $\mathcal{F}_t^{\ti{O},n} \subseteq \mathcal{F}^{\ti O}_t$, which implies $\sigma^2_{\Su} \leq \sigma^2_{\Su,n}$.  

With a finite-memory supplier, our focus is on minimizing the cumulative supply-chain cost \mbox{${\cal C}_{n}(\tpsi) \equiv K_{\I}\,\sigma_{\I}(\tpsi) + K_{\Su}\,\sigma_{\Su, n}(\tpsi)$} within the class of policies $\{\hpsi^m\}$. We obtain the following analogue of \cref{Inventory Approximation}:

\begin{proposition}
\label{Finite Past Cost Approximation} For any $m \in \mathbb{N}$ and for some $C_\kappa > 0$ depending only on $\gamma_\kappa$,
\begin{enumerate}[label={\rm (\roman*)}, leftmargin=3em]
\item $|\sigma^2_{\I}(\tpsi^\infty) - \sigma^2_{\I}(\hpsi^{m})| \le C_\kappa\, m^{-1}$; 
\item  $|\sigma^2_{\Su, n}(\hat{\psi}_m) -  \sigma^2_{\Su*}| \leq C_\kappa\, \Big(\frac{(m + 1)^2 \log n}{n} +\frac{m^2\log (m + \gamma)}{n} + \frac{1}{m}\Big)$;\label{Supplier's Finite Past Error}
\item For a choice of $m(n) \asymp \Big(\frac{\log n}{n}\Big)^{-\frac{1}{3}}$:
\begin{equation}
\label{eq: Main Finite-Past Cost}
    |{\cal C}_n(\hpsi_{m(n)})-{\cal C}^*| \leq C_\kappa\, \Big(\frac{\log n}{n}\Big)^{\frac{1}{3}}.
\end{equation}
\end{enumerate}
 \end{proposition}

Proposition \ref{Finite Past Cost Approximation} underscores the delicate relationship between the supplier's memory $n$ and the complexity $m$ of the retailer's policy $\hpsi^m(z)$. Indeed, from the relation $m(n) \asymp \Big(\frac{\log n}{n}\Big)^{-\frac{1}{3}}$, we see that the retailer should aim to communicate a longer horizon of demand shocks as the supplier widens their forecasting window. However, this dependence is strictly sublinear: if the retailer's order is too complex, the supplier risks overfitting his forecast to a limited order memory, revealing a classical bias-variance tradeoff in designing information delays. Indeed, this will cause the supplier to severely underestimate market volatility from their limited order history and build persistently unreliable forecasts, which, over time, will cause the supplier's inventory costs to skyrocket.

\section{Discussion}

This paper studies how a downstream retailer in a decentralized two-tier supply chain can implicitly transmit demand information to an upstream supplier through the structure of its order stream in the absence of an explicit information-sharing mechanism. We introduce the notion of information delay as a tool for implicit demand sharing, providing an explicit link between optimal implicit information sharing and the group delay of the retailer’s ordering transfer function. We show that pure delay is strictly suboptimal, while fractional-delay mechanisms can reshape the order autocorrelation to improve supplier forecastability and reduce system-wide inventory costs. Using Hardy-space factorization, we develop a tractable family of invertible ARMA policies that approximates the theoretically optimal (but non-rational) limiting filter derived by \cite{caldentey2024information} and preserves its informational delay properties. This construction yields sharp guidance on how policy complexity, as measured by the degrees of the ARMA policies, impacts supply chain costs. We further extend the analysis to memory-constrained suppliers and characterize how the complexity of the retailer’s policy should scale with the supplier’s finite forecasting window, highlighting when, perhaps counterintuitively, increasing policy complexity can become counterproductive.

We close by offering some compelling directions for future research. Our notion of information delay in this paper is characterized by group delay, which can be intuitively seen as the ``first moment'' of a signal's delay structure. While we demonstrate that this delay metric is very well-suited for our analysis here, where costs are expressed as volatilities/variances, it would be interesting to introduce a more comprehensive characterization of the entire delay distribution, which would likely be required to optimize a more complex cost structure. Ongoing work by the authors centers on incorporating more general demand models into our convergence analysis, moving beyond the i.i.d case considered here. Another interesting direction would be to examine if the analysis here can be extended to handle higher-dimensional representations of orders and demand (by perhaps differentiating products). 

\begin{appendices}
\renewcommand{\thefootnote}{\fnsymbol{footnote}}
\setcounter{footnote}{1}
\section{Proofs}
\setcounter{equation}{0}
\renewcommand{\theequation}{A\arabic{equation}}

\noindent {\sc Proof of Theorem \ref{thm:GD_MSFE_Improved}}: Suppose $|Q(z)| \geq c$ for $z \in \mathbb{T}$. We may write:
\begin{equation*}
    Q(z) = \exp\Big(\frac{1}{2\pi} \int_{-\pi}^{\pi} \frac{e^{i\lambda} + z}{e^{i\lambda} - z} \log |Q(e^{i\lambda})|d\lambda\Big) 
\end{equation*}
Now, since $Q(1) = 1$ (bounded inventory condition), we have that:
\begin{align}
    \text{GD}(Q) & = Q'(1) \nonumber  = \frac{Q'(1)}{Q(1)}   = \frac{d}{dz}[\log Q]|_{z = 1}  = \frac{d}{dz}\Big[\frac{1}{2\pi} \int_{-\pi}^{\pi} \frac{e^{i\lambda} + z}{e^{i\lambda} - z} \log |Q(e^{i\lambda})| d\lambda \Big]\Big|_{z = 1} \nonumber \\& = \frac{1}{2\pi} \int_{-\pi}^{\pi} \frac{2 e^{i\lambda}}{(e^{i\lambda} - 1)^2} \log |Q(e^{i\lambda})| d\lambda  = -\frac{1}{2\pi} \int_{-\pi}^{\pi} \Big(1 - \frac{e^{2i\lambda} + 1}{(e^{i\lambda} - 1)^2}\Big) \log |Q(e^{i\lambda})| d\lambda \nonumber \\
    & = -\frac{1}{2}\log \sigma^2_{\Su}(Q) + \frac{1}{2\pi}\int_{-\pi}^{\pi} \frac{\cos \lambda}{1 - \cos \lambda} \log |Q(e^{i\lambda})|d\lambda  \label{eq: Use Kolmogorov}
\end{align}
where \eqref{eq: Use Kolmogorov} follows from Kolmogorov's formula for $\sigma^2_{\Su}(Q)$. 

$\Box$
\vspace{0.5cm}

\noindent {\sc Proof of Corollary \ref{Constant is Controlled}}
We can directly evaluate the integral defining $c$ in \eqref{eq: Constant Definition}:
\begin{equation*}
    c = \frac{1}{2\pi}\int_{-\pi}^{\pi} \frac{\cos \lambda}{1 - \cos \lambda} Q(e^{i\lambda}) d\lambda = \sum_{l = 1}^m \alpha_l\Big(\log \Big|1 - \frac{1}{a_l}\Big| + \text{Re}\Big(\frac{1}{1 - a_l}\Big)\Big) + \sum_{k = 1}^n \log \Big|\sin \frac{\theta_k}{2}\Big| + \frac{1}{2}
\end{equation*}
Clearly the right hand side is jointly decreasing in $\{|a_l|\}_{l = 1}^m$ and blows up at $\alpha_l = 1$ and $\theta_k = 0$ (for any $l \in [m]$ or $k \in [n]$). Since, $\alpha_l \in B(1, \delta)$ by definition, our bound follows. 
$\Box$ \\
\noindent {\sc Proof of Theorem \ref{Inventory_Breakdown}:} By \eqref{eq: Integrability Condition} and Theorem 3.2(v) \cite{fricain2008boundary}, we have that, for all $r \in [0, 1]$, the reproducing kernel:
\begin{equation}
\label{eq: Kernel Definition}
    k_r(z) \equiv \frac{1 - r\tilde{\psi}(r)z\tilde{\psi}(z)}{1 - rz} \in K_{z\tilde{\psi}}
\end{equation}
where $K_{zQ}$ is again the de Branges-Royvnak space corresponding to the function $zQ$ (see \cite{fricain2008boundary} for a definition). Then, from \eqref{eq: Kernel Definition} and Theorem 3.2(v) from \cite{fricain2008boundary}, we have that:
\begin{align*}
    \sigma^2_I(\tilde{\psi}) \langle k_{1}(z), k_1(z) \rangle_{2} & = \lim_{r \to 1} \langle  k_{r}(z), k_r(z) \rangle_{2} \\
    & = \lim_{r \to 1} \langle  k_{r}(z), (I - T_{z\tilde{\psi}}T^{*}_{z\tilde{\psi}})k_r(z) \rangle_{K_{z\tilde{\psi}}}
\end{align*}
For any $\phi \in \mathbb{H}^2$, let $T_{\phi}f = P_{+}(\phi f)$, where $P_{+}$ denotes the orthogonal projection from $L^2(\mathbb{T})$ to $\mathbb{H}^2$, and let $T_{\phi}^{*}$ denote its adjoint. Then, from  \eqref{eq: Kernel Definition} and p.89 of \cite{fricain2008boundary}:
\begin{align*}
    (I - T_{z\tilde{\psi}}T^{*}_{z\tilde{\psi}})k_r(z) & = (I - T_{zQ}T^{*}_{zQ})^2\Big(\frac{1}{1 - rz}\Big) \\
    & = (I - 2T_{zQ}T^{*}_{z\tilde{\psi}} + (T_{z\tilde{\psi}}T^{*}_{z\tilde{\psi}})^2)\Big(\frac{1}{1 - rz}\Big)
\end{align*}
Let $\{b_i\}_{i \in \mathbb{N}}$ denote the Fourier coefficients of the spectral density of $\tilde{\psi}$, i.e:
\begin{equation*}
    |\tilde{\psi}(e^{i\lambda})|^2 = b_0 + \sum_{k \geq 0} b_k(e^{ik\lambda} + e^{-ik\lambda})
\end{equation*}
Hence, we have that:
\begin{align}
    \Big \langle \frac{1}{1-rz}, (T_{z\tilde{\psi}}T^{*}_{z\tilde{\psi}})^2 \Big(\frac{1}{1 - rz}\Big)\Big \rangle & = r^2\tilde{\psi}^2(r) \Big\langle \frac{z\tilde{\psi}(z)}{1-rz}, \frac{z\tilde{\psi}(z)}{1 - rz}\Big\rangle \\
    & = \frac{r^2\tilde{\psi}^2(r)}{1 - r^2} \int_{-\pi}^{\pi}
    \Big(\sum_{i \in \mathbb{Z}} r^i z^i\Big) \Big(\sum_{i \in \mathbb{Z}} b_i z^i \Big) d\lambda \nonumber \\
    & = \frac{r^2\tilde{\psi}^2(r)}{1 - r^2} \sum_{i \in \mathbb{Z}} b_i r^i \equiv \frac{r^2\tilde{\psi}^2(r)}{1 - r^2} \cdot f_Q(r)
\end{align}
where we have defined:
\begin{equation*}
    f_{Q}(re^{i\theta}) = b_0 + 2\sum_{k \geq 1} b_k r^k \cos k\theta
\end{equation*}
with the $f_Q$ denoting the dependence of the spectral density $|\tilde{\psi}(e^{i\lambda})|^2$ only on the outer factor $Q$. Now, observing that:
\begin{equation*}
    T_{z\tilde{\psi}}T^{*}_{z\tilde{\psi}} \Big(\frac{1}{1 - rz}\Big) = \frac{z\tilde{\psi}(z)r\tilde{\psi}(r)}{1 - rz}
\end{equation*}
and putting it all together, we have:
\begin{align*}
     \langle k_{1}(z), k_1(z) \rangle_{2} & = \lim_{r \to 1} \langle  k_{r}(z), (I - T_{z\tilde{\psi}}T^{*}_{z\tilde{\psi}})k_r(z) \rangle_{K_{zQ}} \\
     & = \lim_{r \to 1} \Big\langle  k_{r}(z), \frac{1 - 2z\tilde{\psi}(z)r\tilde{\psi}(r)}{1 - rz} \Big\rangle_{K_{z\tilde{\psi}}} + \frac{r^2\tilde{\psi}^2(r)}{1 - r^2} \cdot f_Q(r)  \\
     & = \lim_{r \to 1} \frac{1 - 2r^2\tilde{\psi}^2(r) +  r^4\tilde{\psi}^4(r)}{1 - r^2} +  \frac{r^2\tilde{\psi}^2(r)(f_Q(r) - r^2\tilde{\psi}^2(r))}{1 - r^2} \\
     & = 0 + \lim_{r \to 1} \frac{r^2\tilde{\psi}^2(r)(f_Q(r) - 1 + 1 - r^2\tilde{\psi}^2(r))}{1 - r^2} \\
     & = 1 + \tilde{\psi}'(1) - \frac{\partial_r f_Q(1)}{2}
\end{align*}
By the bounded inventory condition $\tilde{\psi}(1) = 1$, we have that $\tilde{\psi}'(1) = \frac{\tilde{\psi}'(1)}{\tilde{\psi}(1)} = (\log \tilde{\psi})'(1) = (\log Q)'(1) + (\log B)'(1) + (\log I)'(1) = \text{GD}(Q) + \text{GD}(B) + \text{GD}(I)$ (where the last equality again follows from the bounded inventory condition, applied separately to $Q, B,$ and $I$). $\Box$
\vspace{0.5cm}

\begin{comment}
\iffalse{
\noindent {\sc Proof of Proposition~\ref{prop:info_inner}:}  From Remark~\ref{rem:full_inf}, we have that $\sigma^2_{\Su | \FI}(\tpsi\, I) = |\tpsi(0)\, I(0)|^2$. Since $I(z)$ is analytic in the unit disk  $\mathbb{D}$, the Maximum Modulus Theorem implies that $|I(0)| \leq \max\limits_{z \in \mathbb{T}} |I(z)| = 1$, where the equality follows from the fact that $I(z)$ is inner. The inequality is strict unless I(z)$ is constant in  $\mathbb{D}$. It follows that $
\sigma^2_{\Su | \FI}(\tpsi\, I) \leq |\tpsi(0)|^2= \sigma^2_{\Su | \FI}(\tpsi)$. 
On the other hand, under no information sharing, Kolmogorov's formula in \eqref{eq:Kolmogorov} implies that
\begin{align*}
\sigma^2_{\Su}(\tpsi\,I) &=\exp\left(\frac{1}{2\pi} \int_{-\pi}^\pi \log\left(|\tilde \psi(e^{-i\lambda})\,I(e^{-i\lambda})|^2 \right) \, \D\lambda\right) = \exp\left(\frac{1}{2\pi} \int_{-\pi}^\pi \log\left(|\tilde \psi(e^{-i\lambda})|^2 \right) \, \D\lambda\right) = \sigma^2_{\Su}(\tpsi),
\end{align*}
where the second equality follows from the fact that the inner function $I(z)$ is unimodular on $\mathbb{T}$.
 ~~ $\Box$
\fi}
\vspace{0.5cm}
\end{comment}
%\noindent {\sc Proof of Proposition~\ref{Inventory_Breakdown}:}

\noindent {\sc Proof of \cref{Singular Inner Approximation}:}
To show that the function $\tilde{Q}(z)$ is outer, first note that 
$\tilde{Q}(z)=0$ when $\log(S(z))= - \infty$. Hence, since $S(z) \neq 0$ for $z \in \mathbb{D}$ (as $S$ is singular inner), we likewise have $\tilde{Q}(z) \neq 0$ for $z \in \mathbb{D}$. Moreover, since the singular measure of $S$ is atomic, we have $\log S(z) = \sum_{k \in J} \frac{z + e^{i\lambda_k}}{z - e^{i\lambda_k}}$ for some countable set $J$ --- thus any zeros of $\tilde{Q}$ on $\mathbb{T}$ coincide with roots of $\tilde{Q}$ on $\mathbb{T}$. Finally, since $\log(S(z))=\log(|S(z)|)+i\,\mbox{arg}(z)$, the poles of $\tilde{Q}$ satisfy $\log(|S(z)|)=m>0$.  From the maximum principle applied to the harmonic function $\log |S(z)|$ on $\mathbb{D}$ (recall that $S(z)$ is singular inner and hence $\log S(e^{-i\lambda})$ is purely imaginary, we get that $\log |S(z)| \leq \max_{w \in \mathbb{T}} \lim_{r \to 1^{-}} \log |S(rw)|=0$.  We conclude that $\tilde{Q}$ has no poles and is therefore outer. 

First, we show $\GD(\tilde{Q}) = \phi'(0)$. It is easy to see that:
\begin{align*}
   \GD(\tilde{Q}) & = \tilde{Q}_m'(1)  = m \cdot \frac{S'(z)}{m \cdot S(z)} \Big(1 - \frac{i\log S(z)}{m}\Big)^{m-1}\Big|_{z = 1}  = S'(1)  = \phi'(0)
\end{align*}
where in the penultimate equality, we have used the assumption $\phi(0) = 0$. 

Now, we apply Theorem~\ref{Inventory_Breakdown}. Note, that $\text{GD}(Q\tilde{Q}_m) = \text{GD}(Q) + \text{GD}(\tilde{Q}_m)$ by the bounded inventory condition. Hence, we have that:
\begin{equation*}
    \sigma^2_{\I}(Q\tilde{Q}_m) - \sigma^2_{\I}(Q) = \text{GD}(\tilde{Q}_m) + \frac{\partial_r f_{Q}(1)}{2} - \frac{\partial_r f_{Q\tilde{Q}_m}(1)}{2}. 
\end{equation*}
Now, for any outer transfer function $O \in \mathbb{H}^2$, let $g_O = |O|^2$.  Observe that:
\begin{equation}
\label{eq: Hilbert Representation 1}
    \partial_r f_O(1) = \lim_{r \to 1} \frac{f_{O}(r) - f_{O}(1)}{r - 1} = \frac{d\mathcal{H}[g_O]}{d\lambda}\Big|_{\lambda = 0}
\end{equation}
where $\mathcal{H}[g_O](e^{-i\lambda})$ indicates the Hilbert transform of the spectral density $g_O$ on the unit circle \citep{rosenblum1997hardy}. By the differentiation property of the Hilbert transform, we have that:
\begin{equation}
\label{eq: Hilbert Representation 2}
     \partial_r f_{O}(1) = \mathcal{H}[g_O]'(e^{i0}) = \mathcal{H}[g'_O](e^{i0})
\end{equation}
 Hence, we may write:
\footnotesize
\begin{align}
    |\sigma^2_{\I}(Q\tilde{Q}) - \sigma^2_{\I}(Q) - \text{GD}(\tilde{Q})| & =  \Big|\frac{\partial_r f_{Q}(1)}{2} - \frac{\partial_r f_{Q\tilde{Q}_m}(1)}{2}\Big|  =  \Big|\frac{\mathcal{H}(g_{Q}')(e^{i0}) - \mathcal{H}(g_Q'g_{\tilde{Q}_m})(e^{i0}) - \mathcal{H}(g_Qg'_{\tilde{Q}_m})(e^{i0})}{2} \Big| \label{eq: Hilbert Transform + Product 2} \\
    & \hspace{-2.7cm}= \Big| -\frac{1}{2}\int_{-\pi}^{\pi} g'_Q(e^{i\lambda})\cot \frac{\lambda}{2} d\lambda + \int_{-\pi}^{\pi} g_{\tilde{Q}_m}(e^{i\lambda})g'_Q(e^{i\lambda})\cot \frac{\lambda}{2} d\lambda + \int_{-\pi}^{\pi} g'_{\tilde{Q}_m}(e^{i\lambda})g_Q(e^{i\lambda})\cot \frac{\lambda}{2}d\lambda \Big| \label{eq: Hilbert Def} \\
    & \hspace{-2.7cm} = \Big| - \frac{1}{2}\int_{-\pi}^{\pi}  g'_Q(e^{i\lambda})\Big(1 - \Big(1 + \frac{\phi^2(\lambda)}{m^2}\Big)^{-m}\Big)\cot \frac{\lambda}{2} d\lambda + \frac{1}{m} \int_{-\pi}^{\pi} g_Q(e^{i\lambda}) \Big(1 + \frac{\phi^2(\lambda)}{m^2}\Big)^{-m-1} \phi(\lambda)\phi'(\lambda)\cot \frac{\lambda}{2} d\lambda \Big| \label{eq: Definition of Approximation} \\
    & \hspace{-2.7cm} \leq \frac{1}{2m} \int_{-\pi}^{\pi}  |g'_Q(e^{i\lambda})| \phi^2(\lambda) \cot \frac{\lambda}{2} d\lambda + \frac{1}{m} \int_{-\pi}^{\pi} g_Q(e^{i\lambda}) |\phi(\lambda)\phi'(\lambda)\cot \frac{\lambda}{2}| d\lambda  \label{eq: Mean-Value Theorem 2}
\end{align}
\normalsize 
where \eqref{eq: Hilbert Transform + Product 2} follows from the product rule for differentiation and \eqref{eq: Hilbert Representation 2}, \eqref{eq: Hilbert Def} uses the integral representation of the Hilbert transform on the unit circle, \eqref{eq: Definition of Approximation} is by definition of $\tilde{Q}_m$, and \eqref{eq: Mean-Value Theorem 2} follows from the mean-value theorem (observe $\Big|\frac{d}{dx}(1 + x)^{-m}\Big| \leq m$). Our result then follows by assumption \eqref{eq: Decay Condition on Phase}. $\Box$
\normalsize 

\vspace{0.5cm}

\noindent {\sc Proof of Proposition~\ref{Inventory Approximation}:}
Throughout this proof, we denote:
\begin{align}
    S & = \exp\!\left(\gamma_\kappa\,\frac{z-1}{1 + z}\right) \label{eq: Singular Inner Formula} \\
    \tilde{Q}_m & = \Bigg(1 + \frac{\gamma_\kappa(1 - z)}{m(1 + z)}\Bigg)^{-m}  \label{eq: ARMA Approximation Formula} \\
    Q  & = \frac{1 + z}{2} \label{eq: Overhead Outer}
\end{align}
Now, write $g_O = \|O\|^2$ for the spectral density corresponding to a transfer function $O$. We apply Theorem \ref{Inventory_Breakdown} and observe like in \eqref{eq: Hilbert Representation 1} and \eqref{eq: Hilbert Representation 2}: 
\begin{align}
    |\sigma^2_{\I}(\hat{\psi}_m) - \sigma^2_{\I}(\tilde{\psi}_{\infty})| & = |\text{GD}(\tilde{Q}_m) - \text{GD}(S) + \frac{\partial_r f_{Q}(1)}{2} - \frac{\partial_r f_{Q\tilde{Q}_m}(1)}{2}| \nonumber \\ 
    & = \Big|\frac{\mathcal{H}(g_{Q}')(e^{i0}) - \mathcal{H}(g_Q'g_{\tilde{Q}_m})(e^{i0}) - \mathcal{H}(g_Qg_{\tilde{Q}_m})'(e^{i0})}{2} \Big|\label{eq: Hilbert Transform + Product}  \\
    & = \frac{1}{2}\Big|\int_{-\pi}^{\pi} g'_Q(e^{-i\lambda})\cot \frac{\lambda}{2} d\lambda - \int_{-\pi}^{\pi} g_{\tilde{Q}_m}(e^{i\lambda})g'_Q(e^{-i\lambda})\cot \frac{\lambda}{2} d\lambda - \int_{-\pi}^{\pi} g'_{\tilde{Q}_m}(e^{i\lambda})g_Q(e^{i\lambda})\cot \frac{\lambda}{2}d\lambda \Big| \label{eq: Hilbert Def} \\
    & = \frac{1}{2}\Big|\int_{-\pi}^{\pi} \cos \frac{\lambda}{2}\sin \frac{\lambda}{2}\Big(1 - \Big(1 + \frac{\gamma_{\kappa}^2 \tan^2 \frac{\lambda}{2}}{m^2}\Big)^{-m}\Big)\cot \frac{\lambda}{2} d\lambda  \\
    &  \hspace{2em} -\frac{\gamma_{\kappa}^2}{m}\int_{-\pi}^{\pi} \cos^2 \frac{\lambda}{2}\Big(1 + \frac{\gamma_{\kappa}^2 \tan^2 \frac{\lambda}{2}}{m^2}\Big)^{-m-1}\tan \frac{\lambda}{2} \sec^2 \frac{\lambda}{2}\cot \frac{\lambda}{2} d\lambda \Big| \label{eq: Substitute Expressions} \\
    & \preceq \frac{\gamma_{\kappa}^2}{m} \int_{-\pi}^{\pi} \cos \frac{\lambda}{2}\sin \frac{\lambda}{2} \tan^2 \frac{\lambda}{2} \cot \frac{\lambda}{2} d\lambda + \frac{\gamma^2_{\kappa}}{m} \\
    & \preceq \frac{\gamma^2_{\kappa}}{m}, 
\end{align}
where in \eqref{eq: Hilbert Transform + Product}, we have used that $\text{GD}(\tilde{Q}_m) = \text{GD}(S)$ by Proposition~\ref{Singular Inner Approximation}, applied the Hilbert transform representation of $f_{O}$ (see proof of Proposition~\ref{Singular Inner Approximation}), the commutativity of the Hilbert transform and differentiation, and the product rule, in \eqref{eq: Hilbert Def} we have applied the integral representation of the Hilbert transform on the circle, in \eqref{eq: Substitute Expressions} we have substituted using \eqref{eq: Singular Inner Formula}, \eqref{eq: ARMA Approximation Formula}, and \eqref{eq: Overhead Outer}. Now, it is easy to see that $S(e^{-i\lambda}) = e^{i\tan \frac{\lambda}{2}}$. Hence, we have that:
\begin{equation*}
    \sigma^2_{\I}(\hat{\psi}_m) - \sigma^2_{\I}(\tilde{\psi}_{\infty}) \leq \frac{C_{\kappa}}{m}
\end{equation*}
for some $C_{\kappa} > 0$ depending only on $\kappa$. Moreover, for the part (b) we have:
\begin{align}
    \sigma^2_{\Su}(\hat{\psi}_m) - \sigma^2_{S^{*}} & = \hat{\psi}^2_m(0) - \tilde{\psi}_{\infty}(0) \nonumber \\
    & = \frac{1}{2}\Big((1 - \frac{\gamma_{\kappa}}{m})^{-2m} - \text{exp}(-2\gamma_{\kappa})\Big) \nonumber \\
    & \leq \frac{C_{\kappa}}{m} \label{eq: Infinite-Past Difference}
\end{align}
by a direct application of the mean-value theorem. Putting part \ref{Inventory Error} and \ref{Supplier's Error} together, we obtain the error bound on the cost in part \ref{Total Cost Error}. $\Box$

\vspace{0.5cm}
\noindent {\sc Proof of \cref{Finite Past Cost Approximation}:} Since the error $|\sigma^2_{\I}(\hat{\psi}_m) - \sigma^2_{\I}(\tpsi^{*})|$ is unaffected by the finite-past formulation, we only need to bound the error incurred in the supplier's costs: $|\sigma^2_{\Su, n}(\hat{\psi}_m) -  \sigma^2_{\Su}(\tilde{\psi}_{\infty})|$. We first study $\sigma^2_{\Su, n}(\hat{\psi}_m) -  \sigma^2_{\Su}(\hat{\psi}_m)$. Note, since,
\begin{equation*}
    \max_{m} \int_{-\pi}^{\pi} |\hat{\psi}_m(e^{-i\lambda})|^2 d\lambda < \infty
\end{equation*}
we have that $\max_{n, m} \sigma^2_{\Su, n}(\hat{\psi}_m) < \infty$. For a spectral density $\mu(\theta) = g(\theta)d\theta$, we recall the following expression for $\sigma^2_{\Su, n}(g)$ in terms of Toeplitz determinants:
\begin{equation*}
    \sigma^2_{\Su, n}(g) = \frac{D_{n-1}(g)}{D_{n-2}(g)} 
\end{equation*}
where $D_k(g) = \text{det}\Big([g_{i -  j}]_{i, j =  0}^{k}\Big)$  is the Toeplitz determinant with symbol $g(e^{-i\lambda}) = \sum_{k \in \mathbb{N}} g_k \cos{k\lambda}$ (see e.g. Equation (4.72) in \citealp{gray2004introduction}). Hence, we may write:
\begin{align*}
    \sigma^2_{\Su, n}(\hat{\psi}_m) -  \sigma^2_{\Su}(\hat{\psi}_m)  & \leq |\sigma^2_{\Su, n}(\hat{\psi}_m)| (\log \sigma^2_{\Su, n}(\hat{\psi}_m) - \log \sigma^2_{\Su}(\hat{\psi}_m)) \leq C \Big(\log \frac{D_{n-1}(\hat{g}_m)}{D_{n-2}(\hat{g}_m)} - \frac{1}{\pi} \int \log |\hat{\psi}_m|(e^{-i\lambda})d\lambda \Big)
\end{align*}
where $\hat{g}_m = |\hat{\psi}_m|^2$. Now, applying Theorem 1.1 in \cite{deift2011asymptotics} with:
\begin{align*}
    V(z) & = (z + 1 - \frac{\gamma}{m}(z - 1))^{-m} 
      = \text{exp}\Big(\log (z + 1 - \frac{\gamma}{m}(z - 1))^{-m}\Big)   = \text{exp}\Big(-m\Big(1  + \frac{\gamma}{m} + \sum_{k = 1}^{\infty} \frac{1}{k}\Big(\frac{m - \gamma}{m + \gamma}\Big)^k\Big)\Big)
\end{align*}
we obtain:
\begin{align}
    \sigma^2_{\Su, n}(\hat{\psi}_m) -  \sigma^2_{\Su}(\hat{\psi}_m)  & \leq |\sigma^2_{\Su, n}(\hat{\psi}_m)| (\log \sigma^2_{\Su, n}(\hat{\psi}_m) - \log \sigma^2_{\Su}(\hat{\psi}_m)) \nonumber \\ 
    & \asymp \frac{(m + 1)^2 \log n}{n} \label{eq: Supply Error}
\end{align}
Combining with \eqref{eq: Infinite-Past Difference}, we have:
\begin{equation*}
    |\sigma^2_{\Su, n}(\hat{\psi}_m) -  \sigma^2_{S^{*}}| \preceq \frac{(m + 1)^2 \log n}{n} + \frac{1}{m}
\end{equation*}
which proves the second item in the proposition. Finally combining with the bound on $|\sigma^2_{\I}(\hat{\psi}_m) - \sigma^2_{\I}(\tpsi^{*})|$ in Proposition \ref{Inventory Approximation}, we have:
\begin{equation*}
    |{\cal C}_n(\hpsi_{m(n)})-{\cal C}^*| \leq C_{\kappa}\Big(\frac{(m(n) + 1)^2 \log n}{n} + m(n)^{-1}\Big)
\end{equation*}
Choosing our $m(n)$ to balance terms we obtain the item (iii). $\Box$
\end{appendices}

\bibliographystyle{ormsv080R}
\bibliography{ref}

\end{document}

%% file: Fig_System_OR.tex
\begin{figure}[h]
\begin{center}
\begin{tikzpicture}[scale=1.1]

%%%%
{\iffalse 
\draw (0,0) ellipse (1.2cm and 0.7cm);
\node at (0,0.1) {\scriptsize Manufacturing}; 
\node at (0,-0.3) {\scriptsize Unit};

\draw[<-,>=stealth,line width=1pt] (0.9,0.5) -- (3.4,0.5);
\node at (2.28,0.3) {\etiny manufacturer's Orders}; 
\node at (2.2,0.73) {\footnotesize $O^{\mbox{\etiny S}}_t$}; 

\draw[->,>=stealth,line width=1pt] (0.9,-0.5) -- (3.4,-0.5);
\node at (2.2,-0.27) {\footnotesize $F_t$}; 
\node at (2.28,-0.7) {\etiny Manufacturer's}; 
\node at (2.28,-0.9) {\etiny Shipments}; 
\fi}
%%%%%

\draw [rounded corners] (2.95,-0.7) rectangle (6.2,0.7);
\node at (4.6,0.4) {\bf Supplier};
\node at (4.6,0.05) {\footnotesize (Base-Stock Policy)};

\node at (4.6,-0.43) {\small $S_t=m_t+\sigma_{\ti S}\,\zeta^{\ti S}$};  

\draw[<-,>=stealth,line width=1pt] (6.25,0.5) -- (9.25,0.5);
\node at (7.8,0.3) {\footnotesize Retailer's Orders}; 
\node at (7.8,0.73) {\small $O_t$}; 

\draw[->,>=stealth,line width=1pt] (6.25,-0.5) -- (9.25,-0.5);
\node at (7.8,-0.79) {\footnotesize Supplier's Fulfillment};
\node at (7.88,-1.1) {\ttti{(100\% Service Level)}};

\node at (8,-0.27) {\small $F_t=O_{t-1}$}; 

\draw [rounded corners] (9.3,-0.7) rectangle (12.6,0.7);
\node at (11,0.4) {\bf Retailer}; 
\node at (11,0.05) {\footnotesize (Net Inventory)};
\node at (10.95,-0.43) {\small $I^{\ti R}_t=I^{\ti R}_{t-1}+F_{t}-D_t$};

\draw[<-,>=stealth,line width=1pt] (12.64,0) -- (14.85,0);
\node at (13.8,-0.2) {\footnotesize Demand};
\node at (13.8,0.3) {\small $D_t$}; 
\node at (13.8,-0.5) {\ttti{(Fully Backlogged)}};

\draw (15.8,0) ellipse (0.9cm and 0.7cm);
\node at (15.8,0) {\bf Market}; 

%\draw[dashed] (11,-0.7) -- (11,-1.8);
%\draw[dashed] (5,-1.8) -- (11,-1.8);
%\draw[dashed,->,>=stealth,] (5,-1.8) -- (5,-0.7);
%\node at (8,-1.6) {\footnotesize ${\cal F}^{\ti S}_t$}; 
%\node at (8,-2) {\etiny Information}; 

%\def\shifty{-2.8}

\def\shifty{-1.8}
\def\shiftx{1.5}
\def\shiftD{-0.5}

\node[anchor=west] at (\shiftx,\shifty) {\footnotesize \sf Supplier's information set in period $t$:  ${\cal F}^{\ti S}_t=\sigma\big(O_\tau : \tau \le t\big)$};
\node[anchor=west] at (\shiftx,\shifty+\shiftD) {\footnotesize {\sf Supplier's mean forecast: $m_t=\e[O_{t+1}|{\cal F}^{\ti S}_t]$}};
\node[anchor=west] at (\shiftx,\shifty+2*\shiftD) {\footnotesize \sf Supplier's MSFE: $\sigma_{\ti S}^2=\var[O_{t+1}|{\cal F}^{\ti S}_t]$};

\node[anchor=west] at (\shiftx+7.5,\shifty) {\footnotesize \sf Retailer's inventor cost: $C^{\ti R}=h^{\ti R}\,(I^{\ti R}_t)^+ + b^{\ti R}\,(I^{\ti R}_t)^-$};
\node[anchor=west] at (\shiftx+7.5,\shifty+\shiftD) {\footnotesize \sf Supplier's inventor cost: $C^{\ti S}=h^{\ti S}\,\left(S_{t-1} - O_t\right)^+ + b^{\ti S}\,(O_t - S_{t-1})^+$};
\node[anchor=west] at (\shiftx+7.5,\shifty+2*\shiftD) {\footnotesize \sf Supplier's safety stock factor: $\zeta^{\ti S}=\Phi^{-1}\!\big(\frac{b^{\ti S}}{h^{\ti S}+b^{\ti S}}\big)$};

\end{tikzpicture}
\end{center}
\caption{A two-tier supply chain system.}
 \label{fig:twotierSC}
\end{figure}